\documentclass[12pt]{amsart}
\usepackage[utf8]{inputenc}
\usepackage{euscript}
\usepackage{a4wide}
\usepackage[T1]{fontenc}
\usepackage{lmodern}
\usepackage{amsfonts}
\usepackage{amscd}
\usepackage[matrix,arrow,curve]{xy}
\usepackage{amssymb, amsthm}
\usepackage{amsmath}
\usepackage{pdfpages}
\usepackage{graphicx}
\usepackage{mathtools}
\usepackage{geometry}
\usepackage{tikz}
\usetikzlibrary{calc}
\usepackage{mathrsfs, amssymb,amsthm, amsfonts, amsmath}

\newcounter{cequation}[section]

\theoremstyle{plain}
\newtheorem{thm}[cequation]{Theorem}%[section] %, чтобы нумеровать сначала в каждом разделе
\newtheorem{lemma}[cequation]{Lemma}
\newtheorem{remark}[cequation]{Remark}
\newtheorem{cor}[cequation]{Corollary}
\newtheorem{proposition}[cequation]{Proposition}

\theoremstyle{definition}
\newtheorem{defn}[cequation]{Definition}

\theoremstyle{remark}

\makeatletter\@addtoreset{equation}{section}
\makeatother

\makeatother

\newcommand{\PP}{\mathbb{P}}
\newcommand{\KK}{\mathbb{K}}

\newcommand{\SSS}{\mathfrak{S}}

\newcommand{\mumu}{{\boldsymbol{\mu}}}
\newcommand{\op}{\mathrm{op}}
\newcommand{\sep}{\mathrm{sep}}
\newcommand{\per}{\mathrm{per}}

\newcommand{\Br}{\operatorname{Br}}

\newcommand{\rk}{\operatorname{rk}}

\newcommand{\rkPic}{\rk\Pic}
\newcommand{\WE}{\mathrm{W}(\mathrm{E}_6)}

\newcommand{\Aut}{\operatorname{Aut}}

\newcommand{\Bir}{\operatorname{Bir}}

\newcommand{\Gal}{\operatorname{Gal}}

\newcommand{\Pic}{\operatorname{Pic}}
\newcommand{\Char}{\operatorname{char}}

\usepackage{graphicx}
\graphicspath {{picture/}}
\DeclareGraphicsExtensions{.pdf}

\date{}

\title{Finite Subgroups of Birational Automorphism Groups of Severi--Brauer Surfaces}

\author{Alexandra Sonina}

\address{
Steklov Mathematical Institute of Russian Academy of Sciences, 8 Gubkina St., Moscow 119333, Russia
}

\email{sasha-sonina@mail.ru}

\begin{document}

\begin{abstract}
We prove that a finite group acting birationally on a non-trivial Severi--Brauer surface $S$ over a field $\KK$ of positive characteristic not equal to $3$ is isomorphic either to a subgroup of~$\Aut(S_{\KK^{\per}})$ or to a subgroup of~$\mumu^3_3$, where $\mumu_3$ is a cyclic group of order~$3$ and $\KK^{\per}$ is the perfect closure of a field $\KK$. Combined with results about classification of finite subgroups of group $\Aut(S)$, this gives a classification of finite subgroup of $\Bir(S)$.
\end{abstract}

\maketitle
\tableofcontents

\section{Introduction}

\textit{A Severi--Brauer variety} of dimension $n - 1$ over a field $\mathbb{K}$ is an algebraic variety $X$ which becomes isomorphic to $\mathbb{P}^{n-1}_{\overline{\mathbb{K}}}$ after base change to an algebraic closure $\overline{\mathbb{K}}$ of $\mathbb{K}$. In fact, this is equivalent to requiring it to become trivial over the separable closure $\KK^{\sep}$ of the base field. 
A Severi--Brauer variety is called \textit{non-trivial} if it is not isomorphic to the projective space~$\mathbb{P}^{n-1}_{\mathbb{K}}$ over the base field~$\mathbb{K}$. There is a well-known bijection between Severi--Brauer varieties of dimension~$n - 1$ and central simple algebras of degree~$n$, if one fixes a field~$\mathbb{K}$ and a positive integer~$n$. Furthermore, this bijection preserves automorphism groups, therefore automorphism groups of Severi--Brauer varieties can be described in terms
of the corresponding central simple algebras, see Theorem~E on page~266
of~\cite{Chat}. For more details see e.g.~\cite[~Chapter~5]{G-S}.  

The group of birational transformations of non-trivial Severi--Brauer surfaces has very complicated structure. In the article \cite{BSY} some quotients of $\Bir(S)$ over a perfect fields were described. As a corollary they showed that if there is a point of degree $6$ on $S$, then~$\Bir(S)$ cannot be generated by elements of finite order. The abelianization of $\Bir(S)$ for non-trivial Severi--Brauer surfaces over a perfect fields was found in \cite[Theorem B]{Kurz}. We will approach this group from a different perspective and ask which finite subgroups can act birationally on non-trivial Severi--Brauer surface.

The questions about finite subgroups of $\Bir(S)$ were considered in many papers, for example \cite{Sh20B} and \cite{Vicul}. We will focus on generalizing the results of the article \cite{Sh20B}, which classifies finite subgroups of the birational automorphism groups of surfaces over fields of characteristic zero.

Let $\mumu_n$ denote the cyclic group of order~$n$. 

\begin{thm}[{\cite[Proposition 3.7 and Corollary 4.5]{Sh20B}}]
     \label{Shr3}
   Let~$S$ be a non-trivial Severi--Brauer surface over a field of characteristic zero. Then any finite subgroup of~$\Bir(S)$ is either  conjugate to a subgroup of~$\Aut(S)$ or isomorphic to a subgroup of~$\mumu^3_3$.
\end{thm}

We will prove the version of this theorem for fields of positive characteristic. Firstly, for perfect fields.

\begin{thm}
     \label{main}
   Let~$S$ be a non-trivial Severi--Brauer surface over a perfect field of positive characteristic. Then any finite subgroup of~$\Bir(S)$ is either conjugate to a subgroup of~$\Aut(S)$ or isomorphic to a subgroup of~$\mumu^3_3$.
\end{thm}

\begin{remark}
There is no non-trivial Severi--Brauer surface over a perfect fields of characteristic $3$ (see Proposition \ref{3}).  
\end{remark}

And secondly for arbitrary field of characteristic not equal to $3$. For this result we need the following definition. Let $\KK$ be an arbitrary field, then \textit{the perfect closure} of $\KK$ is a perfect field $\KK^{\per}$ such that~$\KK^{\per}/\KK$ is purely inseparable extension.

\begin{thm}
     \label{main-im}
   Let~$S$ be a non-trivial Severi--Brauer surface over an arbitary field of positive characteristic not equal to 3. Then any finite subgroup of~$\Bir(S)$ is isomorphic either to a subgroup of~$\mumu^3_3$ or to a subgroup of~$\Aut(S_{\KK^{\per}})$, where $S_{\KK^{\per}}$ is still a non-trivial Severi--Brauer surface.
\end{thm}

Let us combine results of the current text with results from the articles \cite{Sh20}, \cite{Sh20B} and~\cite{ASav},~\cite{Son}.

\begin{cor}\label{summerizing}
 Let~$S$ be a non-trivial Severi--Brauer surface over an arbitary field of characteristic not equal to 3 (possibly of characteristic zero). Let $G$ be a finite subgroup of $\Bir(S)$, then either $G$ is a subgroup of $\mumu^3_3$ or a subgroup of $\mumu_3 \times (\mumu_n \rtimes \mumu_3)$, for some integer $n$ and some semidirect product between $\mumu_n$ and $\mumu_3$.  
\end{cor}

Not all groups of the form $\mumu_3 \times (\mumu_n \rtimes \mumu_3)$ can actually act on some non-trivial Severi--Brauer surface and, even more, many of these groups cannot act on a Severi--Brauer surfaces simultaneously, even though each of them may act individually on some different Severi--Brauer surface. For more details see \cite{Son}. 

%Note that \cite[Theorem 1.6]{Son} shows that finite subgroups of both $\Aut(S_{\KK})$ and $\Aut(S_{\KK^{\per}})$ are isomorphic to subgroups of $\mumu_3 \times (\mumu_n \rtimes_{\chi}\mumu_3)$ for various $n$ and $\chi$, where the specific type of semidirect product and the arithmetic conditions on $n$ depend only on the characteristic of the base field. Therefore groups of the forms $\mumu_3 \times (\mumu_n \rtimes_{\chi}\mumu_3)$ are the same for both $\KK$ and $\KK^{\per}$.

%In the articles \cite{ASav} and \cite{Son} in particular all finite supgroups of $\Aut(S)$ were classified. Therefore together with results of this paper one has the classification of all finite subgroups of $\Bir(S)$. 

The following theorem shows, that there are non-trivial Severi--Brauer surfaces on which the group $\mumu^3_3$ acts birationally.

\begin{thm}\label{main-ex}
Let $l$ be a prime number such that $l \neq 3$. There exists a field $\KK$ of characteristic $l$ and a non-trivial Severi--Brauer surface $S$, such that $\mumu_3^3 \subset \Bir(S)$. Moreover, there exists a perfect field of characteristic $l$ with this property as well as an imperfect field of characteristic $l$ with this property.
\end{thm} 

The structure of the paper is as follows: in Section \ref{section:birational}, we remind about birational models of Severi--Brauer surfaces following \cite{Sh20B}; then in Section \ref{section:G-birational}, we recall $G$-birational models of Severi--Brauer surfaces again following \cite{Sh20B}, but adapting the proofs to perfect fields of arbitrary characteristic rather than fields of characteristic $0$. In Section~\ref{section:bounds}, Theorem~\ref{main} is proven. In Section~\ref{imperfect}, Theorem~\ref{main-im} and Corollary \ref{summerizing} are proven and it is shown that in the case when the base field has characteristic $3$ more work need to be done. Finally, in Section~\ref{example} we construct examples which prove Theorem \ref{main-ex}.

\smallskip
\textbf{Acknowledgements.}
I would like to express my sincere gratitude to Constantin Shramov for proposing the problem, for numerous indispensable discussions for his careful reading. I am thankful to Andrey Trepalin for several helpful discussions and explanations about del Pezzo surfaces. I am also thankful to Elias Kurz for the useful reference on the current understanding of the group of birational automorphisms of non-trivial Severi--Brauer surfaces.

This work was performed at the Steklov International Mathematical Center and partially supported by the Ministry of Science and Higher Education of the Russian Federation (agreement no. 075-15-2025-303) and by the Theoretical Physics and Mathematics Advancement Foundation «BASIS».

\section{Birational models of Severi--Brauer surfaces}
\label{section:birational}

In this section, we recall the results of \cite{Sh20B} concerning surfaces that may be birationally equivalent to a non-trivial Severi--Brauer surface.

%The following general result is sometimes referred to as the theorem of Lang and Nishimura.

\begin{thm}[{see e.g.~\cite[Proposition~IV.6.2]{Kollar-RatCurves}}]
\label{theorem:Lang-Nishimura}
Let $X$ and $Y$ be smooth projective varieties over an arbitrary field $\KK$.
Suppose that $X$ is birational to $Y$. Then $X$ has a~$\KK$-point if and only if $Y$ has a $\KK$-point.
\end{thm}

\begin{cor}[{\cite[Corollary 2.2]{Sh20B}}]\label{corollary:Lang-Nishimura}
Let $X$ and $Y$ be smooth projective varieties over an arbitrary field $\KK$, and let $r$ be a positive integer.
Suppose that~$X$ is birational to $Y$, and~$X$ has a point of degree not divisible by~$r$.
Then $Y$ has a point of degree not divisible by~$r$.
\end{cor}

The following result concerning Severi--Brauer surfaces is well-known.

\begin{thm}[{see e.g.~\cite[Theorem~53(2)]{Kollar-SB}}]
\label{theorem:SB-point-degree}
Let $S$ be a non-trivial Severi--Brauer surface over an arbitrary field. Then $S$ does not contain points of degree $d$ not divisible by $3$.
\end{thm}

Recall the definition of a del Pezzo surface and a conic bundle.

\begin{defn}
A del Pezzo surface is a projective smooth surface with an ample anticanonical class.
\end{defn}

For a del Pezzo surface $S$, by its degree we mean its (anti)canonical degree~$K_S^2$.

\begin{defn}
A conic bundle is a proper flat morphism $\pi: S\to C$ of smooth varieties 
such that it is of relative dimension one and the anti-canonical divisor $-K_S$ is relatively ample.
\end{defn}

In our cases the dimension of $S$ is always $2$, so $C$ is always a smooth curve. Note that every fiber will be a conic, since $-K_S$ is relatively ample. 

\begin{cor}[{\cite[Corollary 2.4]{Sh20B}}]\label{corollary:many-dP}
Let $S$ be a non-trivial Severi--Brauer surface over an arbitrary field. Then
$S$ is not birational to any conic bundle, and not birational to any del Pezzo surface of degree
different from $3$, $6$, and $9$.
\end{cor}

\begin{cor}[{\cite[Corollary 2.5]{Sh20B}}]\label{corollary:dP-large-Picard-rank}
Let $S$ be a non-trivial Severi--Brauer surface over an arbitrary field. Then~$S$
is not birational to any del Pezzo surface of degree $6$ of Picard rank greater than~$2$,
and not birational to any del Pezzo surface of degree $3$ of Picard rank greater than~$3$.
\end{cor}

To proceed we will need the following general fact about non-trivial Severi--Brauer surfaces. First, recall definitions.

\begin{defn}
Let $d\le 6$ and $P_1,\ldots,P_d$ are distinct $\KK^{\sep}$-points on $\PP^2_{{\KK^{\sep}}}$ over an separably closed field ${\KK^{\sep}}$. We will say that they are \emph{in general position} if no three of them are contained in a line, and all six (in the case $d=6$) are not contained in a conic~(cf.~\mbox{\cite[Remark~IV.4.4]{Manin}}).
\end{defn}

\begin{defn}
Let~$P$ be a separable point of degree $d$ on a Severi--Brauer
surface $S$ over a field~$\KK$, we will say that~$P$ is in general position if the $d$ points
of the set~$P_{\KK^{\sep}}\subset\PP^2_{\KK^{\sep}}$ are in general position. 
\end{defn}

\begin{lemma}[{\cite[Lemma 2.8]{Sh20B}}]\label{lemma:SB-general-position}
Let $S$ be a non-trivial Severi--Brauer surface over a perfect field~$\KK$,
and let~\mbox{$P\in S$} be a point of degree $d$. Suppose that $d=3$ or $d=6$. Then~$P$ is a point in general position.
\end{lemma}

Let $S$ be a non-trivial Severi--Brauer surface over a perfect field~$\KK$, and let~$A$
be the corresponding central simple algebra.
We denote by $S^{\op}$ the Severi--Brauer surface corresponding to the algebra opposite to~$A$.
There are two classes of interesting birational maps from $S$ to $S^{\op}$ which we describe below
(cf. \cite[Lemma~3.1]{IskovskikhTregub}, or  cases~(c) and~(e) of~\cite[Theorem~2.6(ii)]{Is-UMN}).

Let $P$ be a point of degree $3$ on $S$. Then $P$ is a point in general position by Lemma~\ref{lemma:SB-general-position}.
Blowing up $P$ and blowing down the proper transforms of the three lines on~${S_{{\KK^{\sep}}}\cong\PP^2_{{\KK^{\sep}}}}$
passing through the pairs of points of $P_{{\KK^{\sep}}}$,
we obtain a birational map $\tau_P$ to another Severi--Brauer surface $S'$.
This map is given by the linear system of conics passing through~$P$.

Similarly, let $P$ be a point of degree $6$ on $S$. Then $P$ is a point in general position by Lemma~\ref{lemma:SB-general-position}.
Blowing up $P$ and blowing down the proper transforms of the six conics on~\mbox{$S_{{\KK^{\sep}}}\cong\PP^2_{{\KK^{\sep}}}$}
passing through the quintuples of points of $P_{{\KK^{\sep}}}$,
we obtain a birational map $\eta_P$ to another Severi--Brauer surface $S'$.
This map is given by the linear system of quintics singular at~$P$~(or, in other words, at each point of~$P_{{\KK^{\sep}}}$).

In both of the above cases one has $S'\cong S^{\op}$. This follows from a well-known general
fact about the degrees of birational maps between Severi--Brauer varieties with given classes,
see for instance~\cite[Exercise~3.3.7(iii)]{GSh}.
For more details on the maps $\tau_P$, see~\mbox{\cite[\S2]{Corn}}.

\begin{remark}\label{remark:Amitsur}
It is known that $S$ and $S^{\op}$ are the only Severi--Brauer surfaces birational to~$S$, see for instance~\cite[Exercise~3.3.6(v)]{GSh}.
\end{remark}

\begin{cor}[{\cite[Corollary 2.11]{Sh20B}}]\label{corollary:Weinstein}
Let $S$ be a non-trivial Severi--Brauer surface over a perfect field.
Then~$S$ is not birational to any del Pezzo surface $S'$ of degree~$3$ or~$6$
with~\mbox{$\rkPic(S')=1$}.
\end{cor}

The part of Corollary~\ref{corollary:Weinstein} concerning del Pezzo surfaces of degree $3$ also follows
from~\mbox{\cite[Chapter~V]{Manin}}.
The part concerning del Pezzo surfaces of degree $6$ can be obtained from
\cite[\S2]{IskovskikhTregub}.

Corollaries~\ref{corollary:dP-large-Picard-rank} and \ref{corollary:Weinstein} show that a del Pezzo surface
of degree $6$ birational to a non-trivial Severi--Brauer
surface must have Picard rank equal to~$2$, and a del Pezzo surface
of degree $3$ birational to a non-trivial Severi--Brauer
surface must have Picard rank equal to~$2$ or~$3$.
In the next section we will obtain further restrictions on such surfaces provided
that they are $G$-minimal with respect to some finite group~$G$.

\section{$G$-birational models of Severi--Brauer surfaces}
\label{section:G-birational}

In this section, following \cite{Sh20B}, we study finite groups acting on surfaces that are birationally equivalent to a non-trivial Severi--Brauer surface.

We start with del Pezzo surfaces of degree $6$ over a field $\KK$. Recall that over $\KK^{\sep}$ a del Pezzo surface of degree $6$ is unique up to isomorphism, and its automorphism group is isomorphic to $((\KK^{\sep})^*)^2\rtimes (\SSS_3\times\mumu_2)$.
More details on this can be found in~\cite[~Theorem~8.4.2]{Dolgachev} and \cite{Zaitsev}.

Given a del Pezzo surface $S'$ of degree $6$ over a perfect field $\KK$, we will call $\SSS_3\times\mumu_2$ its \emph{Weyl group}.
For every element $\theta\in\Aut(S')$ we will refer to the image of $\theta$
under the composition of the embedding $\Aut(S')\hookrightarrow\Aut(S'_{{\KK^{\sep}}})$ with the natural homomorphism
$$
\Aut(S'_{{\KK^{\sep}}})\to\SSS_3\times\mumu_2
$$
as the image of $\theta$ in the Weyl group. Similarly, we will consider the image of the Galois group
$\Gal(\KK^{\sep}/\KK)$ in the Weyl group.

Let $\sigma\colon\PP^2\dasharrow\PP^2$ be the \emph{standard Cremona involution},
that is, a birational involution acting as
$$
[x:y:z]\mapsto\left[\frac{1}{x}:\frac{1}{y}:\frac{1}{z}\right]
$$
in some homogeneous coordinates $x$, $y$, and $z$. This involution becomes regular
on the del Pezzo surface of degree $6$ obtained as a blow up of $\PP^2$ at the points
$[1:0:0]$, $[0:1:0]$, and $[0:0:1]$. Let $\hat{\sigma}$ be its image in the Weyl group $\SSS_3\times\mumu_2$.
Then $\hat{\sigma}$ is the generator of the center $\mumu_2$ of the Weyl group.
If one thinks about~\mbox{$\SSS_3\times\mumu_2$}
as the group of symmetries of a regular hexagon, then $\hat{\sigma}$ is an involution
that interchanges the opposite sides of the hexagon or, in other words, a rotation by~$180^\circ$.
If $S'$ is a del Pezzo surface of degree $6$ over some perfect field and $\theta$ is its automorphism,
we will say that $\theta$ is \emph{of Cremona type} if its image in the Weyl group
coincides with $\hat{\sigma}$.

Let us provide several statements from \cite{Sh20B} with their original proofs. The proofs remain valid over any perfect field, so the assumption that the base field has characteristic zero was unnecessary.

\begin{lemma}[{\cite[Lemma 3.1]{Sh20B}}]\label{lemma:Cremona-type}
Let $S'$ be a del Pezzo surface of degree $6$ over a perfect field $\KK$ of characteristic $l$, and let $\theta$
be its automorphism of Cremona type. Then $\theta$ is an involution that has exactly four fixed points on~$S'_{{\KK^{\sep}}}$ if $l \neq 2$ and exactly one if $l = 2$.
\end{lemma}

\begin{proof}
Using a birational contraction $S'_{{\KK^{\sep}}}\to\PP^2_{{\KK^{\sep}}}$,
consider the automorphism $\theta$
as a birational automorphism of $\PP^2_{{\KK^{\sep}}}$.
Then $\theta$ can be represented as a composition of the standard Cremona
involution $\sigma$ with some element of the standard torus acting on $\PP^2_{{\KK^{\sep}}}$.
Thus, one can choose homogeneous coordinates $x$, $y$, and $z$
on $\PP^2_{{\KK^{\sep}}}$ so that $\theta$ acts as
$$
[x:y:z]\mapsto\left[\frac{\alpha}{x}:\frac{\beta}{y}:\frac{\gamma}{z}\right]
$$
for some non-zero $\alpha,\beta,\gamma\in{\KK^{\sep}}$.
Therefore, $\theta$ is conjugate to the involution~$\sigma$ via the automorphism
$$
[x:y:z]\mapsto [\sqrt{\alpha}x:\sqrt{\beta}y:\sqrt{\gamma}z].
$$
This shows that $\theta$ has the same number of fixed points on $\PP^2_{{\KK^{\sep}}}$ as $\sigma$,
while it is easy to see that if $l \neq 2$ the fixed points of the latter are the four points
$$
[1:1:1],\quad [-1:1:1], \quad [1:-1:1],\quad [1:1:-1].
$$
Note that if $l = 2$, then the only fixed point is $[1:1:1]$.
It remains to notice that a fixed point of $\theta$ on $S'_{{\KK^{\sep}}}$ cannot be contained
in a $(-1)$-curve, and thus all of them are mapped to fixed points of $\theta$ on~$\PP^2_{{\KK^{\sep}}}$.
\end{proof}

\begin{lemma}[{\cite[Lemma 3.2]{Sh20B}}]\label{lemma:dP6-rk-1}
Let $S$ be a non-trivial Severi--Brauer surface over a perfect field $\KK$, and let
$S'$ be a del Pezzo surface of degree $6$ over $\KK$. Suppose that there exists a finite group
$G\subset\Aut(S')$ such that $\rkPic(S')^G=1$. Then $S'$ is not birational to~$S$.
\end{lemma}

\begin{proof}
Suppose that $S'$ is birational to $S$.
We know from Corollary~\ref{corollary:dP-large-Picard-rank} that~\mbox{$\rkPic(S')\le 2$}.
Hence by Corollary~\ref{corollary:Weinstein} one has~\mbox{$\rkPic(S')=2$}.
Since $S'$ does not have a structure of a conic bundle by Corollary~\ref{corollary:many-dP}, we see that
$S'$ has a contraction on a del Pezzo surface of larger degree.
Again by Corollary~\ref{corollary:many-dP},
this means that $S'$ is a blow up of a Severi--Brauer surface at a point of degree~$3$.
Hence the image $\Gamma$ of the Galois group~\mbox{$\Gal({\KK^{\sep}}/\KK)$} in the Weyl group~\mbox{$\SSS_3\times\mumu_2$}
contains an element of order~$3$.

Since the cone of effective curves
on $S'$ has two extremal rays and $\rkPic(S')^G=1$, we see that $G$ must contain an element
whose image in the Weyl group has order~$2$. On the other hand,
the image of $G$ in the Weyl group commutes with~$\Gamma$. Since $\hat{\sigma}$
is the only element of order $2$ in~\mbox{$\SSS_3\times\mumu_2$}  that commutes with an element
of order~$3$, we conclude that~$G$
must contain an element $\theta$
of Cremona type. The element~$\theta$ has exactly $4$ or~$1$ fixed points on $S'_{{\KK^{\sep}}}$ by Lemma~\ref{lemma:Cremona-type}.
This implies that there is a point of degree not divisible by~$3$ on $S'$.
Thus the assertion follows from Corollary~\ref{corollary:Lang-Nishimura}  and Theorem~\ref{theorem:SB-point-degree}.
\end{proof}

Now we deal with del Pezzo surfaces of degree $3$, i.e. smooth cubic surfaces in $\PP^3$. Recall that for a del Pezzo surface
$S'$ of degree $3$ over a field $\KK$ the action of the groups~\mbox{$\Aut(S')$} and $\Gal({\KK^{\sep}}/\KK)$ on $(-1)$-curves on $S'$ defines homomorphisms
of these groups to the Weyl group~$\WE$. Furthermore, the homomorphism
$\Aut(S')\to\WE$ is an embedding. The order of the Weyl group~$\WE$ equals~$2^7\cdot 3^4\cdot 5$.
We refer the reader to~\cite[Corollary~8.2.40]{Dolgachev}
and~\cite[Chapter~IV]{Manin} for details.

Now we will give a proof of the following lemma, which is not just a modification of similar statment from the article \cite{Sh20B}.

\begin{lemma}\label{lemma:dP3-not-3-group}
Let $S$ be a non-trivial Severi--Brauer surface over a perfect field $\KK$ of positive characteristic $l$, and let $S'$ be a del Pezzo surface of degree $3$ over $\KK$. Suppose that there exists a non-trivial automorphism $g\in\Aut(S')$ such that the order of $g$ is not a power of~$3$. Then~$S'$ is not birational to $S$.
\end{lemma}

\begin{proof}
We may assume that the order $p=\mathrm{ord}(g)$ is prime. Thus one has $p=2$ or~\mbox{$p=5$}, since the order of the Weyl group $\WE$ is not divisible by primes greater than~$5$.

 The action of $g$ on $S'_{{\KK^{\sep}}}$ is of one of the three types listed in \cite[Table 7]{Dolgachev-Duncan}. We may assume that $S'_{{\KK^{\sep}}}$ is in a form as in \cite[Table 8]{Dolgachev-Duncan}.
 
Firstly, suppose $g$ has order $5$. Then $g$ has type $5A$ and $l \neq 5$. Then $S'_{{\KK^{\sep}}}$ is the Clebsch cubic. Let us find a fixed point locus. Since $l \neq 5$, then there are exactly $4$ fixed points~${[1:\xi:\xi^2:\xi^3:\xi^4]}$, where $\xi$ is any primitive root of unity of degree $5$. 

Secondly, suppose $g$ has order $2$. Then $g$ can have either type $2A$ or $2B$. 

Suppose $l \neq 2$. If $g$ has type $2A$, then the fixed point locus consists of an elliptic curve and an isolated point. If $g$ has type $2B$ then the fixed point locus consists of a $(-1)$-curve and three isolated points. 

Suppose $l = 2$. If $g$ has type $2B$, then the fixed point locus consists of a $(-1)$-curve.

If $g$ has type $2A$, then the form from \cite[Table 8]{Dolgachev-Duncan} imposes some conditions on the cubic. However, by \cite[Theorem 9.2]{Dolgachev-Duncan} the center of the involution of type $2A$ is Eckardt point and moreover by \cite[Proposition 9.1]{Dolgachev-Duncan} there are three $g$-invariant lines passing through the center of the involution.

Note that a $\Gal({\KK^{\sep}}/\KK)$-invariant $(-1)$-curve on $S'$ always contains a $\KK$-point, since such a curve is a line
in the anticanonical embedding of $S'$. Note also that if there are exactly three invariant lines passing through a single point, then this point  must be fixed under action of $\Gal({\KK^{\sep}}/\KK)$.

Therefore, in each of the above three cases we find a $\Gal({\KK^{\sep}}/\KK)$-invariant set of points on $S'_{\bar{\KK}}$ of cardinality coprime to~$3$. Thus the assertion follows from Corollary~\ref{corollary:Lang-Nishimura}  and Theorem~\ref{theorem:SB-point-degree}.

\end{proof}

\begin{cor}[{\cite[Corollary 3.4]{Sh20B}}]\label{corollary:dP3-rkPic-3}
Let $S$ be a non-trivial Severi--Brauer surface over a perfect field $\KK$ of positive characteristic, and let
$S'$ be a del Pezzo surface of degree~$3$ over $\KK$ birational to $S$.
Suppose that there exists a subgroup $G\subset\Aut(S')$ such that~${\rkPic(S')^G=1}$. Then~$\rkPic(S')=3$.
\end{cor}

\begin{proof}
We know from Corollaries~\ref{corollary:dP-large-Picard-rank} and
\ref{corollary:Weinstein} that
either $\rkPic(S')=2$, or $\rkPic(S')=3$.
Suppose that $\rkPic(S')=2$, so that the cone of effective curves
on $S'$ has two extremal rays. Since $\rkPic(S')^G=1$, the group $G$ must contain an element of even order.
Therefore, the assertion follows from Lemma~\ref{lemma:dP3-not-3-group}.
\end{proof}

\begin{cor}[{\cite[Corollary 3.5]{Sh20B}}]\label{corollary:dP3-3-group}
Let $S$ be a non-trivial Severi--Brauer surface over a perfect field $\KK$ of positive characteristic, and let $S'$ be a del Pezzo surface of degree $3$ over~$\KK$ birational to~$S$.
Let~\mbox{$G\subset\Aut(S')$} be a subgroup such that $\rkPic(S')^G=1$. Then~$G$ is isomorphic to
a subgroup of~$\mumu_3^3$.
\end{cor}

\begin{proof}
We know from Corollary~\ref{corollary:dP3-rkPic-3} that $\rkPic(S')=3$.
Hence $S'$ is a blow up of a Severi--Brauer surface at two points of degree $3$.
This means that the image $\Gamma$ of the Galois group~\mbox{$\Gal({\KK^{\sep}}/\KK)$} in
the Weyl group $\WE$ contains an element conjugate to
$$
\gamma=(123)(456)
$$
in the notation of~\cite[\S4]{Trepalin-cubic}. 
We may assume that $\Gamma$ contains $\gamma$ itself, so that
the image of~$G$ in $\WE$ is contained in the centralizer $Z(\gamma)$
of $\gamma$. By~\mbox{\cite[Proposition~4.5]{Trepalin-cubic}}
one has
$$
Z(\gamma)\cong(\mumu_3^2\rtimes\mumu_2)\times \SSS_3.
$$
Note that since both the notation and the result about the centeralizer are formulated for the abstract group$\WE$, then the characteristic of the base field doesn't matter.  
On the other hand, we know from Lemma~\ref{lemma:dP3-not-3-group} that the order of $G$ is a power of $3$.
Since the Sylow $3$-subgroup of $Z(\gamma)$ is isomorphic to $\mumu_3^3$, the required
assertion follows.
\end{proof}

Let us summarize the results of this section.

\begin{proposition}[{\cite[Proposition 3.7]{Sh20B}}]\label{proposition:summary}
Let $S$ be a non-trivial Severi--Brauer surface over a perfect field $\KK$ of positive characteristic, and let
$G\subset\Bir(S)$ be a finite subgroup. Then $G$ is conjugate either to a subgroup of
$\Aut(S)$, or to a subgroup of $\Aut(S^{\op})$, or to a subgroup of
$\Aut(S')$, where $S'$ is a del Pezzo surface of degree $3$ over $\KK$ birational to~$S$
such that~\mbox{$\rkPic(S')=3$} and~\mbox{$\rkPic(S')^G=1$}.
In the latter case $G$ is isomorphic to
a subgroup of~$\mumu_3^3$.
\end{proposition}

\begin{proof}
Regularizing the action of $G$ and running a $G$-Minimal Model Program (see~\mbox{\cite[Theorem~1G]{Iskovskikh80}}),
we obtain a $G$-surface $S'$ birational to $S$, such that~$S'$ is either a del Pezzo surface with
$\rkPic(S')^G=1$, or a conic bundle. The case of a conic bundle is impossible by Corollary~\ref{corollary:many-dP}.
Thus by Corollary~\ref{corollary:many-dP} and Lemma~\ref{lemma:dP6-rk-1}
we conclude that $S'$ is a del Pezzo surface of degree $9$ or $3$.
In the former case $S'$ is a Severi--Brauer surface itself, so
that $S'$ is isomorphic either to $S$ or to $S^{\op}$ by Remark~\ref{remark:Amitsur}.
In the latter case
we have $\rkPic(S')=3$ by Corollary~\ref{corollary:dP3-rkPic-3}. Furthermore, in this case~$G$ is isomorphic to
a subgroup of~$\mumu_3^3$ by
Corollary~\ref{corollary:dP3-3-group}.
\end{proof}

\section{Proof of Theorem \ref{main}}
\label{section:bounds}

In this section we prove Theorem~\ref{main}.
We start by recalling a result about $G$-invariant point.

\begin{proposition}[{\cite[Proposition A.5]{Son}}]
\label{corollary:SB-Aut}
Let $S$ be a non-trivial Severi--Brauer surface over a field of characteristic not equal to $3$, and
let $G\subset\Aut(S)$ be a finite subgroup. Then there exists a $G$-invariant separable point in general position of degree $3$ on $S$.
\end{proposition}

Recall from~\S\ref{section:birational} that a  Severi--Brauer surface $S$ is birational to the surface~$S^{\op}$. In particular, the groups $\Bir(S)$ and $\Bir(S^{\op})$ are (non-canonically) isomorphic, and the group $\Aut(S^{\op})$ is (non-canonically) realized as a subgroup of~\mbox{$\Bir(S)$}.

Proposition~\ref{corollary:SB-Aut} has the following geometric consequence analogous to the one established in \cite{Sh20B} for surfaces over fields of characteristic zero.

\begin{cor}[{\cite[Corollary 4.5]{Sh20B}}]\label{corollary:conjugate}
Let $S$ be a non-trivial Severi--Brauer surface over a perfect field, and
let $G\subset\Aut(S^{\op})$ be a finite subgroup. Then $G$ is conjugate to a subgroup of~$\Aut(S)$.
\end{cor}

\begin{proof}
By Proposition~\ref{corollary:SB-Aut} the group $G$ has
an invariant separable point $P$ in general position of degree $3$ on $S^{\op}$.
Blowing up~$P$ and blowing down the proper transforms of the three lines on $S^{\op}_{{\KK^{\sep}}}\cong\PP^2_{{\KK^{\sep}}}$
passing through the pairs of the three points of $P_{{\KK^{\sep}}}$,
we obtain a (regular) action of $G$ on the surface~$S$ together with a $G$-equivariant
birational map~${\tau_P\colon S^{\op}\dasharrow S}$, cf.~\S\ref{section:birational}.
This means that $G$ is conjugate to a subgroup of~$\Aut(S)$.
\end{proof}

Now we prove Theorem~\ref{main}.

\begin{proof}[Proof of Theorem~\ref{main}]
We know from Proposition~\ref{proposition:summary} that $G$ is either conjugate to a subgroup of
$\Aut(S')$, where $S'\cong S$ or $S'\cong S^{\op}$ or isomorphic to subgroup of $\mumu^3_3$. In the first and the last cases we are done.
In the second case~$G$ is conjugate to a subgroup of~$\Aut(S)$ by Corollary~\ref{corollary:conjugate}.
\end{proof}

\section{Imperfect fields and fields of characteristic $3$}
\label{imperfect}

In this section we generalize our result to the case of imperfect fields of characteristic not equal to~$3$.

\begin{defn}
Let $\KK$ be an arbitrary field, then \textit{the perfect closure} of $\KK$ is a perfect field~$\KK^{\per}$ such that~$\KK^{\per}/\KK$ is purely inseparable extension.
\end{defn}

The existence and uniqueness of the perfect closure is guaranteed by \cite[tag 046W]{Stack}.

\begin{lemma}\label{perfectization}
Let $l$ be a prime number such that $l \neq 3$. Let $\KK$ be a field of characteristic~$l$. Let $S$ be a non-trivial Severi--Brauer surface. Then $S_{\KK^{\per}}$ is also a non-trivial Severi--Brauer surface.
\end{lemma}

\begin{proof}
Suppose that $S_{\KK^{\per}}$ is a trivial Severi--Brauer surface. Then it has a point, which was defined over some finite purely inseparable extension of $\KK$. All finite purely inseparable extension of $\KK$ has degree equal to some power of $l$, so the point has degree equal to some power of $l$. Since $\Char(\KK) \neq 3$, this gives a contradiction with Theorem \ref{theorem:SB-point-degree}. Therefore~$S_{\KK^{\per}}$ is a non-trivial Severi--Brauer surface.
\end{proof}

Now we are ready to prove Theorem \ref{main-im}.

\begin{proof}[Proof of Theorem \ref{main-im}]
Let $\KK^{\per}$ be the perfect closure of the field $\KK$. By Lemma \ref{perfectization} one has that $S_{\KK^{\per}}$ is non-trivial. Now we get a finite group~$G$ acting on $S_{\KK^{\per}}$, i.e. on a non-trivial Severi--Brauer variety over a perfect field. Then the result follows from Theorem~\ref{main}. 
\end{proof}

\begin{proof}[Proof of Corollary \ref{summerizing}]
If $\Char(\KK) = 0$ the result follows by combining results of articles~\cite{Sh20B} and \cite{Sh20}. Suppose $\Char(\KK) > 0$ and not equal to $3$. Then by Theorem \ref{main-im} one has that $G$ is isomorphic either to a subgroup of $\mumu^3_3$ or to a subgroup of $\Aut(S_{\KK^{per}})$. If~$G$ is isomorphic to a subgroup of $\Aut(S_{\KK^{per}})$ then by~\cite[~Theorem 1.4]{ASav}, there exists an integer~$n$ and a semidirect product between $\mumu_n$ and $\mumu_3$ such that $G$ is isomorphic to a subgroup of~${\mumu_3 \times (\mumu_n \rtimes \mumu_3)}$. Actually, semidirect product has a specific type and there are arithmetic conditions on $n$ depend only on the characteristic of the base field. For more details see \cite{Son}.  
\end{proof}

Finally, let us notice that non-trivial Severi--Brauer surfaces do not exist over perfect fields of characteristic $3$. 

%therefore there is no point in asking the question in this case.

\begin{proposition}\label{3}
Let $S$ be a Severi--Brauer surface over a perfect field $\KK$ of characteristic~$3$. Then $S$ is trivial.
\end{proposition} 

\begin{proof}
Recall that $\Omega^1_{\KK}$ denote the module of relative differentials, see e.g. \cite[Chapter~II,~Paragraph 8]{Har}.

%i.e. a vector space over $\KK$, generated by $dk$, where $k\in \KK$ and factorized by the following relations: 
%\begin{enumerate}
%\item [1] $d(k + k') - dk - dk'$, where $k,k' \in \KK$
%\item [2] $d(kk') - kdk' - k'dk$, where $k,k' \in \KK$
%\item [3] $d1$, where $1$ is a neutral element. 
%\end{enumerate}

If $\KK$ is a perfect field of characteristic $p$, then for every element $k' \in \KK$ there exists an element $k \in \KK$, such that $k' = k^p$. Therefore $$dk' = pk^{p-1}dk = 0,$$ since the characteristic of the field $\KK$ is $p$. It means that $\Omega^1_{\KK}$ is trivial if $\KK$ is a perfect field of positive characteristic.

By \cite[Theorem 9.2.4]{G-S} one has a surjective map from $\Omega^1_{\KK}$ on to the $3$-torsion of Brauer group $\Br(\KK)$ of the base field. Note that class of the surface $S$ lies in the $3$-torsion of~$\Br(\KK)$. Since~$\Omega^1_{\KK}$ is trivial, the $3$-torsion of~$\Br(\KK)$ is also trivial. Therefore $S$ is trivial.
\end{proof}

Since for imperfect field $\KK$ of characteristic $3$ every Severi--Brauer surface $S$ becomes trivial over $\KK^{\per}$, there is not much we can say about $\Bir(S)$ in this case.

\section{Examples}
\label{example}
In this section we provide a construction similar to the main construction from \cite[~Paragraph 3]{ShCub} that proves 
Theorem~\ref{main-ex}.

Let $\Bbbk$ be an algebraically closed field of characteristic $l$, and let  $\KK=\Bbbk(\lambda,\mu)$, where~$\lambda$ and $\mu$ are transcendental variables.
Consider the cubic surface $S$ over $\KK$ given in 
the projective space $\PP^3$ with homogeneous coordinates 
$x$, $y$, $z$, and $t$ by the equation
\begin{equation*}
\lambda x^3+\lambda^2 y^3+\mu z^3+\mu^2 t^3=0.
\end{equation*}
It is straightforward to check that $S$ has no $\KK$-points. 
Since $l \neq 3$ there exists a non-trivial cubic root of unity in $\Bbbk$. Denote it by $\omega$.
Consider the linear transformations 
\begin{equation*}
\begin{aligned}
&\sigma_1\colon [x:y:z:t]\mapsto [\omega x:y:z:t],\\
&\sigma_2\colon [x:y:z:t]\mapsto [x:\omega y:z:t],\\
&\sigma_3\colon [x:y:z:t]\mapsto [x:y:\omega z:t].
\end{aligned}
\end{equation*}
Then $\sigma_1$, $\sigma_2$, and $\sigma_3$ generate a subgroup $G$ of 
$\Aut(S)$ isomorphic to $\mumu_3^3$.

Recall that the action of the groups~\mbox{$\Aut(S)$} 
and $\Gal\left({\KK^{\sep}}/\KK\right)$ 
on $(-1)$-curves on $S_{\KK^{\sep}}$ defines homomorphisms
of these groups to the Weyl group~$\WE$. Furthermore, the homomorphism
$\Aut(S)\to\WE$ is an embedding, see~\cite[Corollary~8.2.40]{Dolgachev}.
Thus we can identify the group $\Aut(S)$ with its image in $\WE$. 
In the notation of \cite{Carter}, 
the elements~$\sigma_i$, for~$1\leqslant i\leqslant 3$, 
have type~$\mathrm{A}_2^3$; this follows from the 
fact that the quotient of $S_{{\KK^{\sep}}}$ by the group 
$G_i\cong \mumu_3$ generated 
by the element $\sigma_i$ is isomorphic to~$\PP^2$,
and thus~\mbox{$\rkPic(S_{{\KK^{sep}}})^{G_i}=1$}.

Let $\gamma_{\lambda}$ and $\gamma_{\mu}$ be the images in $\WE$ 
of generators of the Galois 
groups~\mbox{$\Gal\left(\KK\left(\sqrt[3]{\lambda}\right)/\KK\right)$} 
and 
$\Gal\left(\KK\left(\sqrt[3]{\mu}\right)/\KK\right)$.
Then $\gamma_{\lambda}$ and $\gamma_{\mu}$ are elements of order $3$ that 
commute with each other. 
Furthermore, 
the $27$ lines on $S_{{\KK^{sep}}}$ are given by equations
$$
\begin{aligned}
\sqrt[3]{\lambda}x+\omega^i\sqrt[3]{\lambda^2}y=
\sqrt[3]{\mu}z+\omega^j\sqrt[3]{\mu^2}t=0,& \\
\sqrt[3]{\lambda}x+\omega^i\sqrt[3]{\mu}z=
\sqrt[3]{\lambda^2}y+\omega^j\sqrt[3]{\mu^2}t=0,& \\
\sqrt[3]{\lambda}x+\omega^i\sqrt[3]{\mu^2}t=
\sqrt[3]{\lambda^2}y+\omega^j\sqrt[3]{\mu}z=0,& \\ 
i,j\in\{0,1,2\}.&
\end{aligned}
$$
Thus each of the~$27$ lines on $S_{{\KK^{\sep}}}$
is defined over the field $\KK\left(\sqrt[3]{\lambda},\sqrt[3]{\mu}\right)$. 
Hence  the image of the whole Galois group $\Gal\left({\KK^{\sep}}/\KK\right)$ 
in the Weyl group $\WE$ is the group $\Gamma\cong\mumu_3^2$ 
generated by~$\gamma_{\lambda}$ and~$\gamma_{\mu}$. 

Denote the lines in the first set by $L^1_{i,j}$, the lines in the second set by  $L^2_{i,j}$, and the lines in the third set by $L^3_{i,j}$. Under the action of group $\Gamma$ there is one orbit consisting of $9$ lines $L^1_{i,j}$ and $6$ orbits each consisting of $3$ lines $L^n_{i,j}$ with $i+j$ is some constant modulo $3$ and $n = 2,3$.  

Note that $L^n_{i,j} \cap L^n_{k,l} = \varnothing$ if and only if $i \neq k$ and $j \neq l$ for $1 \leqslant n \leqslant 3$. Therefore lines in any orbit of length $3$ do not intersect each others. Note that $L^2_{i,j} \cap L^3_{k,l} \neq \varnothing$ if and only if $i+j \equiv k+l \pmod 3$. Therefore the six lines in the union of two $\Gamma$-orbits of length $3$, namely $$L^2_{0,0}, L^2_{1,2}, L^2_{2,1}, L^3_{1,0},L^3_{2,2}, L^3_{0,1}$$ do not intersect each other.

We can contract these two $\Gamma$-orbits over $\KK$. After the contraction we get a del Pezzo surface of degree $9$. This del Pezzo surface has no $\KK$-points by Theorem \ref{theorem:Lang-Nishimura}. Therefore we contracted $S$ onto a non-trivial Severi--Brauer variety $S'$. Note that $\mumu_3^3$ acts birationally on $S'$, therefore~$S'$ is an example we needed.  

The field $\KK$ used in the construction above is imperfect. To obtain an example over a perfect field it suffices to consider the base change $S'_{\KK^{\per}}$. The group $\mumu_3^3$ still acts birationally on $S'_{\KK^{\per}}$, and $S'_{\KK^{\per}}$ is non-trivial, since $\Char(\KK) \neq 3$.

\end{document}